\documentclass[12pt]{amsart}

\usepackage{amssymb}

\def\cal{\mathcal}
\def\opn#1#2{\def#1{\operatorname{#2}}} 
\opn\chara{char}
\opn\length{\ell}
\opn\pd{pd}
\opn\rk{rk}
\opn\projdim{proj\,dim}
\opn\injdim{inj\,dim}
\opn\rank{rank}
\opn\depth{depth}
\opn\grade{grade}
\opn\height{height}
\opn\embdim{emb\,dim}
\opn\codim{codim}

\opn\Tr{Tr}
\opn\bigrank{big\,rank}
\opn\superheight{superheight}\opn\lcm{lcm}
\opn\trdeg{tr\,deg}%
\opn\reg{reg}
\opn\lreg{lreg}
\opn\skel{skel}
\opn\com{com}
\opn\div{div}
\opn\Div{Div}
\opn\cl{cl}
\opn\Cl{Cl}
\opn\Spec{Spec}
\opn\Supp{Supp}
\opn\supp{supp}
\opn\Sing{Sing}
\opn\Ass{Ass}
\opn\Ann{Ann}
\opn\Rad{Rad}
\opn\Soc{Soc}
\opn\Ker{Ker}
\opn\Coker{Coker}
\opn\Im{Im}
\opn\Hom{Hom}
\opn\Tor{Tor}
\opn\Ext{Ext}
\opn\End{End}
\opn\Aut{Aut}
\opn\id{id}

\opn\nat{nat}
\opn\pff{proof}
\opn\Pf{Pf}
\opn\GL{GL}
\opn\SL{SL}
\opn\mod{mod}
\opn\ord{ord}
\opn\aff{aff}
\opn\con{conv}
\opn\relint{relint}
\opn\st{st}
\opn\lk{lk}
\opn\cn{cn}
\opn\core{core}
\opn\vol{vol}
\opn\link{link}
\opn\star{star}
\opn\gr{gr}

\def\pot#1#2{#1[\kern-0.28ex[#2]\kern-0.28ex]}

\opn\dirlim{\underrightarrow{\lim}}
\opn\inivlim{\underleftarrow{\lim}}
\def\Implies{\ifmmode\Longrightarrow \else
     \unskip${}\Longrightarrow{}$\ignorespaces\fi}
\def\implies{\ifmmode\Rightarrow \else
     \unskip${}\Rightarrow{}$\ignorespaces\fi}
\def\iff{\ifmmode\Longleftrightarrow \else
     \unskip${}\Longleftrightarrow{}$\ignorespaces\fi}

\let\:=\colon
\newtheorem{Theorem}{Theorem}[section]
\newtheorem{Lemma}[Theorem]{Lemma}
\newtheorem{Corollary}[Theorem]{Corollary}

\newtheorem{Remark}[Theorem]{Remark}

\let\epsilon\varepsilon
\let\phi=\varphi
\let\kappa=\varkappa
\opn\initial{in}
\opn\inim{inm}
\opn\rev{rev}
\opn\Gin{Gin}
\opn\Lex{Lex}
\opn\Shift{Shift}
\opn\shift{shift}
\opn\rate{rate}
\opn\Mon{Mon}
\opn\lex{lex}
\opn\rev{rev}
\opn\red{red}
\opn\max{max}
\opn\min{min}
\opn\initial{in}
\opn\Ker{Ker}
\opn\GL{GL}
\opn\proj{proj}
\begin{document}
\title{
Unmixed bipartite graphs
and
sublattices of
the Boolean lattices}

\author{J\"urgen Herzog}
\address{
J\"urgen Herzog,
Fachbereich Mathematik und Informatik, Universit\"at
Duisburg--Essen, Campus Essen, 45117 Essen, GERMANY
}
\email{juergen.herzog@uni-essen.de}

\author{Takayuki Hibi}
\address{
Takayuki Hibi,
Department of Pure and Applied Mathematics, Graduate School of Information Science and Technology,
Osaka University, Toyonaka, Osaka 560-0043,
JAPAN
}
\email{hibi@math.sci.osaka-u.ac.jp}

\author{Hidefumi Ohsugi}
\address{
Hidefumi Ohsugi,
Department of Mathematics, College of Science, 
Rikkyo University,
Tokyo 171-8501, JAPAN}
\email{ohsugi@rkmath.rikkyo.ac.jp}

\begin{abstract}
The correspondence between 
unmixed bipartite graphs and 
sublattices of the Boolean lattice is discussed.
By using this correspondence, we show the existence of
squarefree quadratic initial ideals of toric ideals
arising from minimal vertex covers of unmixed bipartite graphs.
\end{abstract}

\maketitle

\section*{Introduction}

Let $G$ be a finite graph on the vertex set $[N]=\{1,\ldots,N\}$
with no loops, no multiple edges and no isolated vertices.
A {\it vertex cover} of $G$ is a subset $C \subset [N]$ such that,
for each edge $\{i,j\}$ of $G$, one has either $i \in C$ or $j \in C$.
Such a vertex cover $C$ is called {\it minimal} if
no subset $C' \subsetneqq C$
is a vertex cover of $G$.
We say that a finite graph $G$ is {\it unmixed} if all minimal vertex covers of $G$
have the same cardinality.
Let $A = K[z_1, \ldots, z_N]$ the polynomial ring in $N$ variables
over a field $K$.
The {\em edge ideal} of $G$ is
the monomial ideal $I(G)$ of $A$ generated by those quadratic
monomials $z_iz_j$ such that $\{ i, j \}$ is an edge of $G$.
It is well-known that
the primary decomposition of the edge ideal of $G$ is
$$
I(G) = \bigcap_{C \in {\cal M}(G)} \left< z_i \ | \ i \in C \right>
$$
where ${\cal M}(G)$ is the set of all minimal vertex covers of $G$.
We say that $G$ is {\em Cohen--Macaulay} (over $K$) if
the quotient ring $A / I(G)$ is Cohen--Macaulay.
Every Cohen--Macaulay graph is unmixed.
A graph theoretical characterization of
Cohen--Macaulay bipartite graphs was given in \cite{HH}
and that of unmixed bipartite graphs was given in \cite{V}.

In Sections 1 and 2, we study the correspondence between
unmixed bipartite graphs and sublattices of the
Boolean lattice ${\cal L}_n$ on $\{x_1,\ldots,x_n\}$:
\begin{itemize}
\item
There exists a one-to-one correspondence between
unmixed bipartite graphs on $\{x_1,\ldots,x_n\} \cup \{y_1,\ldots,y_n\}$
and sublattices ${\cal L}$ of ${\cal L}_n$
with $\emptyset \in {\cal L}$ and $\{x_1,\ldots,x_n\} \in {\cal L}$.
(Theorem \ref{unmixed}.)
\item
There exists a one-to-one correspondence between
Cohen--Macaulay bipartite graphs on $\{x_1,\ldots,x_n\} \cup \{y_1,\ldots,y_n\}$
and full sublattices of ${\cal L}_n$.
(Theorem \ref{CM}.)
\end{itemize}

In Section 3, we study toric ideals arising from 
the set of minimal vertex covers of unmixed bipartite graphs.
Let $G$ be an unmixed bipartite graph on the vertex set
$\{x_1,\ldots,x_n\} \cup \{y_1,\ldots,y_n\}$
and let
$K[x_1,\ldots,x_n,y_1,\ldots,y_n]$ the polynomial ring in $2n$ variables
over a field $K$ with each $\deg x_i = \deg y_i=1$.
We associate each minimal vertex cover $C$ of $G$ with the squarefree monomial
$
u_C = \prod_{v  \in C} v
\in K[x_1,\ldots,x_n,y_1,\ldots,y_n]$ of degree $n$.
Let ${\cal R}_G$ denote the semigroup ring
generated by all monomials $u_C$ with $C \in {\cal M}(G)$
over $K$.
Let $S_G = K[\{ z_C \}_{C \in {\cal M}(G)}]$ denote the polynomial ring 
of $|{\cal M}(G)|$ variables over $K$.
The {\it toric ideal} $I_G$ of ${\cal R}_G$ is the kernel of the
surjective homomorphism $\pi: S_G \rightarrow  {\cal R}_G$
defined by $\pi (z_C) = u_C$.
In Section 3, 
by using the correspondence given in Section 1,
we show that:
\begin{itemize}
\item
The toric ideals arising from an unmixed bipartite graph possesses
a squarefree quadratic initial ideal.
(Theorem \ref{quadratic}.)
\end{itemize}

\section{Minimal vertex covers of unmixed bipartite graphs}

First we recall a fact stated in \cite[p.300]{HH}.
Let $G$ be a bipartite graph without isolated vertices
and let $V(G) = \{x_1,\ldots,x_m\} \cup \{y_1,\ldots,y_n\}$ denote the set of vertices of $G$.
Suppose that $G$ is unmixed.
Since $\{x_1,\ldots,x_m\}$ and $\{y_1,\ldots,y_n\}$ are minimal vertex covers
of $G$, we have $m = n$.
Moreover, by virtue of the ``marriage theorem," 
we may assume that $\{x_i,y_i\} \in E(G)$ for all $i$.

Thanks to this fact, it follows that
each minimal vertex cover of $G$ is of the form
$\{x_{i_1},\ldots,x_{i_s},y_{i_{s+1}},\ldots,y_{i_n}\}$
where 
$\{i_1,\ldots, i_n \} = [n] $.
For a minimal vertex cover $C=\{x_{i_1},\ldots,x_{i_s},y_{i_{s+1}},\ldots,y_{i_n}\}$
of $G$, we set $\overline{C} = \{x_{i_1},\ldots,x_{i_s}\}$.
Let ${\cal L}_n$ denote the Boolean lattice on the set $\{x_{1},\ldots,x_{n}\}$
and let 
$${\cal L}_G = \{\overline{C} \ | \ C \mbox{ is a minimal vertex cover of  } G\}
\ (\subset {\cal L}_n).$$

\begin{Remark}
{\rm
Let $G$ and $G'$ be unmixed bipartite graphs on
$\{x_1,\ldots,x_n\} \cup \{y_1,\ldots,y_n\}$.
\begin{enumerate}
\item[(i)]
Both
$\emptyset$ and $\{x_{1},\ldots,x_{n}\}$
belong to ${\cal L}_G$.
\item[(ii)]
If $G \neq G'$, then we have $I(G) \neq I(G')$.
Hence ${\cal L}_G \neq {\cal L}_{G'}$ follows from
the primary decomposition of the edge ideals.
\end{enumerate}
}
\end{Remark}

\begin{Theorem}
\label{unmixed}
Let ${\cal L}$ be a subset of ${\cal L}_n$.
Then
there exists a (unique) unmixed bipartite graph $G$
on $\{x_1,\ldots,x_n\} \cup \{y_1,\ldots,y_n\}$
such that ${\cal L} = {\cal L}_G$
if and only if
$\emptyset$ and
$\{x_{1},\ldots,x_{n}\}$ belong to ${\cal L}$ and
${\cal L}$ is a sublattice of ${\cal L}_n$.
\end{Theorem}

\begin{proof}
{\bf (``Only if")}
Suppose that both
$C= \{x_{i_1},\ldots,x_{i_s},y_{i_{s+1}},\ldots,y_{i_n}\}$
and
$C'=\{x_{j_1},\ldots,x_{j_t},y_{j_{t+1}},\ldots,y_{j_n}\}$
are minimal vertex covers of $G$.
Then
$$\{ y_k \ | \ x_k \notin \overline{C} \cap \overline{C'} \} = \{y_{i_{s+1}},\ldots,y_{i_n}\} \cup \{ y_{j_{t+1}},\ldots,y_{j_n}\},$$
$$\{ y_k \ | \ x_k \notin \overline{C} \cup \overline{C'} \}=\{y_{i_{s+1}},\ldots,y_{i_n}\} \cap \{ y_{j_{t+1}},\ldots,y_{j_n}\}.$$

First we show that $\overline{C} \cap  \overline{C'} \in {\cal L}_G$, that is,
$C_1=( \overline{C} \cap \overline{C'} ) \cup \{ y_k \ | \ x_k \notin \overline{C} \cap \overline{C'} \}$
is a minimal vertex cover of $G$.
Suppose that an edge $\{x_i,y_j\}$ of $G$ 
satisfies
$y_j \notin \{ y_k \ | \ x_k \notin \overline{C} \cap \overline{C'} \}= \{y_{i_{s+1}},\ldots,y_{i_n}\} \cup \{ y_{j_{t+1}},\ldots,y_{j_n}\}$.
Since $C$ (resp. $C'$) is a vertex cover of $G$, we have $x_i \in \overline{C}$ (resp. $x_i \in 
\overline{C'}$).
Hence $x_i \in \overline{C} \cap \overline{C'}$.
Thus $C_1$ is a minimal vertex cover of $G$.

Second we show that $\overline{C} \cup  \overline{C'} \in {\cal L}_G$, that is,
$C_2= (\overline{C} \cup \overline{C'})
\cup \{ y_k \ | \ x_k \notin \overline{C} \cup \overline{C'} \}$
is a minimal vertex cover of $G$.
Suppose that an edge $\{x_i,y_j\}$ of $G$ 
satisfies
$x_i \notin \overline{C} \cup \overline{C'}$.
Since $C$ (resp. $C'$) is a vertex cover of $G$, we have $y_j \in \{y_{i_{s+1}},\ldots,y_{i_n}\}$ (resp. $y_j \in 
\{ y_{j_{t+1}},\ldots,y_{j_n}\}$).
Thus
$
y_j \in \{y_{i_{s+1}},\ldots,y_{i_n}\} \cap \{ y_{j_{t+1}},\ldots,y_{j_n}\}
= \{ y_k \ | \ x_k \notin \overline{C} \cup \overline{C'} \}$
and hence $C_2$ is a minimal vertex cover of $G$.

\bigskip

{\bf (``If")}
For each element $S \in {\cal L}$, let $S^*$ denote
the set $\{y_j \ | \ x_j \notin S\}$.
Let $I$ be an ideal
$
\bigcap_{S \in {\cal L}} \left< S \cup S^* \right>
$.
We will show that
there exists an unmixed bipartite graph $G$ such that
$I = \left< x_i y_j \ | \ \{x_i,y_j\} \in E(G) \right>$.

Since $\emptyset \in {\cal L}$ and $\{x_{1},\ldots,x_{n}\} \in {\cal L}$,
$I \subset \left< x_i y_j \ | \ 1 \leq i,j \leq n \right> $.
Suppose that a monomial $M$ of degree $\geq 3$ belongs to  the minimal set of generators of $I$.

If $M = x_i x_j u$ where $i \neq j$ and $u$ is a (squarefree) monomial,
then there exist $S, S' \in {\cal L}$ such that
$x_i \in S \setminus S'$,
$x_j \in S' \setminus S$,
$u \notin \left< S \cup S^* \right>$ and
$u \notin \left< S' \cup {S'}^* \right>$.
Since ${\cal L}$ is a sublattice of ${\cal L}_n$,
$S \cap S' \in {\cal L}$.
Note that $(S \cap S')^* = S^* \cup {S'}^*$.
Hence we have
$$
I \subset
\left<
(S \cap S') \cup (S^* \cup {S'}^*)
\right>
.$$
However,
none of the variables in $M$ appears in
the set
$(S \cap S') \cup (S^* \cup {S'}^*).$

If $M = y_i y_j u$ where $i \neq j$ and $u$ is a (squarefree) monomial,
then there exist $S, S' \in {\cal L}$ such that
$y_i \in S^* \setminus {S'}^*$,
$y_j \in {S'}^* \setminus S^*$,
$u \notin \left< S \cup S^* \right>$ and
$u \notin \left< S' \cup {S'}^* \right>$.
Since ${\cal L}$ is a sublattice of ${\cal L}_n$,
$S \cup S' \in {\cal L}$.
Note that $(S \cup S')^* = S^* \cap {S'}^*$.
Hence we have
$$
I \subset
\left<
(S \cup S') \cup (S^* \cap {S'}^*)
\right>
.$$
However,
none of the variables in $M$ appears in
the set
$(S \cup S') \cup (S^* \cap {S'}^*).$

Thus the minimal set of generators of $I$ is a subset of
$\{x_i y_j \ | \ 1 \leq i,j \leq n\}$
and hence there exists a bipartite graph $G$ such that
$I = I(G)$.
Since the primary decomposition of the edge ideal $I(G)$ of $G$ is
$I = \bigcap_{C \in {\cal M}(G)} \left< C \right>$,
it follows that ${\cal M}(G) = \{ S \cup S^* \ | \ S \in {\cal L} \}$.
Thus we have
${\cal L} = {\cal L}_G$.
Since the cardinality of
each $S \cup S^*$ with $S \in {\cal L}$
is $n$, $G$ is unmixed as desired.
\end{proof}

\section{Minimal vertex covers of Cohen--Macaulay bipartite graphs}

Let, as before, $G$ be a finite graph on $[N]$ and
$A = K[z_1, \ldots, z_N]$ the polynomial ring in $N$ variables
over a field $K$.  The {\em edge ideal} of $G$ is
the monomial ideal $I(G)$ of $A$ generated by those quadratic
monomials $z_iz_j$ such that $\{ i, j \}$ is an edge of $G$.
We say that $G$ is {\em Cohen--Macaulay} (over $K$) if
the quotient ring $A / I(G)$ is Cohen--Macaulay.
Every Cohen--Macaulay graph is unmixed.

Given a finite poset (partially ordered set) 
$P = \{ p_1, \ldots, p_n \}$,
we introduce the bipartite graph $G_P$
on the vertex set 
$\{ x_1, \ldots, x_n \} \cup \{ y_1, \ldots, y_n \}$
whose edges are those 2-element subsets $\{ x_i, y_j \}$
with $p_i \leq p_j$.  
In particular for each $i \in [n]$ the edge $\{ x_i, y_i \}$
belongs to $G_P$.
It is known \cite{HH} that $G_P$
is Cohen--Macaulay.  Conversely, given a Cohen--Macaulay
bipartite graph $G$, there is a finite poset $P$ with
$G = G_P$.

A subset $\alpha \subset P$ is called a {\em poset ideal}
of $P$ if $\alpha$ enjoys the property that
if $p_i \in \alpha$ and $p_j \leq p_i$, then
$p_j \in \alpha$.  In particular
the empty set and $P$ itself
are poset ideals of $P$.
For each poset ideal $\alpha$ of $P$,
we set
$\alpha_x =
\{ x_i \, | \, p_i \in \alpha \}$
and
$\alpha_y =
\{ y_j \, | \, p_j \not\in \alpha \}$.

\begin{Lemma}
\label{Oberwolfach}
The set $\alpha_x \cup \alpha_y$ is a minimal vertex cover
of $G_P$.  Conversely, every minimal vertex cover of $G_P$
is of the form $\beta_x \cup \beta_y$
for some poset ideal $\beta$ of $P$.
\end{Lemma}
 
\begin{proof}
Let $\alpha$ be a poset ideal of $P$.  We show that
$C = \alpha_x \cup \alpha_y$ is a minimal vertex cover
of $G_P$.  Let $\{ x_i, y_j \}$ be an edge of $G$.
Then $p_i \leq p_j$.  Suppose $x_i \not\in \alpha_x$.  
Then $p_i \not\in \alpha$.  Since $\alpha$ is a poset ideal of $P$,
it follows that $p_j \not\in \alpha$.  Thus $y_j \in \alpha_y$.
Hence $C$ is a vertex cover of $G$.  Since
$G_P$ is unmixed and $|C| = n$, it follows that
$C$ is a minimal vertex cover. 

Conversely, given a minimal vertex cover
$C = \{ x_{i_1}, \ldots, x_{i_s} \} 
\cup \{ y_{i_{s+1}}, \ldots, y_{i_{n}} \}$ of $G_P$,
where $\{ i_1, \ldots, i_n \} = [n]$, we prove that
$\alpha = \{ p_{i_1}, \ldots, p_{i_s} \}$ 
is a poset ideal of $P$.
Let $p_{i_j} \in \alpha$ and $p_{a} < p_{i_j}$ in $P$.
Then $\{ x_a, y_{i_j} \}$ is an edge of $G_P$.
Suppose $p_a \not\in \alpha$.  Then $x_a \not\in C$.
Since $x_{i_j} \in C$, one has $y_{i_j} \not\in C$.   
Thus neither $x_a$ nor $y_{i_j}$ belongs to $C$.
However, $\{ x_a, y_{i_j} \}$ is an edge of $G_P$.
Thus $C$ cannot be a vertex cover of $G_P$.
Hence $p_a \in \alpha$.  Consequently, 
$\alpha$ is a poset ideal of $G_P$, as desired.
\end{proof}

Let, as before, ${\mathcal L}_n$ denote the Boolean lattice
on $\{x_1, \dots, x_n\}$.  A sublattice ${\mathcal L}$
of ${\mathcal L}_n$ is called {\em full} if the rank of
${\mathcal L}$ is equal to $n$.  Here the rank of 
${\mathcal L}$ is defined to be the nonnegative integer
$\ell - 1$, where
$\ell$ is the maximal cardinality of chains 
(totally ordered subsets) of ${\mathcal L}$. 

Let $P$ be a finite poset with $|P| = n$ and
${\mathcal J}(P)$ the set of all poset ideals of $P$.
It turns out that the subset ${\mathcal J}(P)$ 
of ${\mathcal L}_n$ is a full
sublattice of ${\mathcal L}_n$.
Conversely, the classical fundamental structure theorem 
for finite distributive lattices 
\cite[pp. 118--119]{Hibibook}
guarantees that every full sublattice of ${\mathcal L}_n$
is of the form ${\mathcal J}(P)$ for a unique poset
$P$ with $|P| = n$.

\begin{Theorem}
\label{CM}
A subset ${\mathcal L}$ of ${\mathcal L}_n$
is a full sublattice of ${\mathcal L}_n$
if and only if
there exists a Cohen--Macaulay bipartite graph $G$
on $\{x_1,\ldots,x_n\} \cup \{y_1,\ldots,y_n\}$
with ${\mathcal L} = {\mathcal L}_G$.
\end{Theorem}

\begin{proof}
{\bf (``If")}
Let $G$ be a Cohen--Macaulay bipartite graph
on the set $\{x_1,\ldots,x_n\} \cup \{y_1,\ldots,y_n\}$ and 
$P$ a poset with $G = G_P$, where $|P| = n$. 
Lemma \ref{Oberwolfach} says that
${\mathcal L}_G$ coincides with ${\mathcal J}(P)$.
Thus ${\mathcal L}_G$ is a full sublattice of ${\mathcal L}_n$.

{\bf (``Only if")}
Suppose that ${\mathcal L}$ is a full sublattice of ${\mathcal L}_n$.
One has ${\mathcal L} = {\mathcal J}(P)$ for a unique poset
$P$ with $|P| = n$.  Let $G = G_P$.
Then $G$ is a Cohen--Macaulay bipartite graph.  
Lemma \ref{Oberwolfach} says that
${\mathcal L}_G$ coincides with ${\mathcal J}(P)$.
Thus ${\mathcal L}_G = {\mathcal L}$, as required.
\end{proof}

\section{Toric ideals arising from minimal vertex covers}

Let $G$ be an unmixed bipartite graph on the vertex set
$\{x_1,\ldots,x_n\} \cup \{y_1,\ldots,y_n\}$
and let
$K[x_1,\ldots,x_n,y_1,\ldots,y_n]$ the polynomial ring in $2n$ variables
over a field $K$ with each $\deg x_i = \deg y_i=1$.
Let ${\cal M}(G)$ denote the set of all minimal vertex covers of $G$.
We associate each minimal vertex cover $C$ of $G$ with the squarefree monomial
$
u_C = \prod_{v  \in C} v
\in K[x_1,\ldots,x_n,y_1,\ldots,y_n]$ of degree $n$.
Let ${\cal R}_G$ denote the semigroup ring
generated by all monomials $u_C$ with $C \in {\cal M}(G)$
over $K$.
Let $S_G = K[\{ z_C \}_{C \in {\cal M}(G)}]$ denote the polynomial ring 
of $|{\cal M}(G)|$ variables over $K$.
The {\it toric ideal} $I_G$ of ${\cal R}_G$ is the kernel of the
surjective homomorphism $\pi: S_G \rightarrow  {\cal R}_G$
defined by $\pi (z_C) = u_C$.

\begin{Theorem}
\label{quadratic}
Let $G$ be an unmixed bipartite graph.
Then the toric ideal $I_G$ of ${\cal R}_G$ has a squarefree quadratic initial ideal
with respect to a reverse lexicographic order.
\end{Theorem}

\begin{proof}
Let $G_0$ denote the (unmixed) bipartite graph with the edge set
$E(G) =  \{ \{x_1,y_1 \} ,\ldots, \{x_n,y_n \}\}$.
Then ${\cal L}_{G_0}$ is the Boolean lattice on $\{x_1,\ldots,x_n\}$.
It is known \cite{Hibipaper} that 
the reduced Gr\"obner basis of toric ideal of ${\cal R}_{G_0}$ 
with respect to a suitable reverse lexicographic order is
$$
{\cal G}_0=
\{
\ \underline{z_C z_{C'}} - z_{C \cap C'} z_{C \cup C'} \ | \ C, C' \in {\cal M}(G) , \ C \neq C' \ 
\}
$$
where the initial monomial of each binomial of ${\cal G}_0$ is 
the first monomial.

Let $G$ be an unmixed bipartite graph on the vertex set
$\{x_1,\ldots,x_n\} \cup \{y_1,\ldots,y_n\}$.
Then ${\cal L}_G$ is a sublattice of ${\cal L}_{G_0}$.
Hence we have the following:
\begin{itemize}
\item[(i)]
$I_G = I_{G_0} \cap S_G$
(by \cite[Proposition 4.13 (a)]{Stu});
\item[(ii)]
If $C$ and $C'$ belong to ${\cal M}(G)$,
then ${C \cap C'}$ and ${C \cup C'}$ belong to ${\cal M}(G)$.
Thus, if $z_C z_{C'} \in S_G$, then we have $z_{C \cap C'} z_{C \cup C'} \in S_G$.
\end{itemize}
Thanks to the elimination property above, 
${\cal G}_0 \cap S_G$ 
is a Gr\"obner basis of the toric ideal $I_G$ of ${\cal R}_G$ as desired.
\end{proof}

\begin{Corollary}
Let $G$ be an unmixed bipartite graph.
Then the semigroup ring ${\cal R}_G$ is normal and Koszul.
\end{Corollary}

\end{document}